# OPTIMAL HOEFFDING BOUNDS FOR DISCRETE REVERSIBLE MARKOV CHAINS

By Carlos A. León[1] and François Perron[2]

*Universidad de Concepción and Université de Montréal*

We build optimal exponential bounds for the probabilities of large deviations of sums $\sum_{k=1}^{n} f(X_k)$ where $(X_k)$ is a finite reversible Markov chain and $f$ is an arbitrary bounded function. These bounds depend only on the stationary mean $\mathbb{E}_\pi f$, the end-points of the support of $f$, the sample size $n$ and the second largest eigenvalue $\lambda$ of the transition matrix.

**0. Introduction.** Consider an ergodic Markov chain $(X_k)$ with finite state space $E$, transition matrix $P$ and stationary distribution $\pi$. Let $f : E \to \mathbb{R}$ satisfy $\min f(E) = 0$, $\max f(E) = 1$ and let $\mu = \int f \, d\pi$. From the weak law of large numbers we know that the empirical mean $n^{-1} S_n = n^{-1} \sum_{k=1}^{n} f(X_k)$ converges to $\mu$ in probability. This result is the working principle behind all Markov chain Monte Carlo (MCMC) integration techniques. The basis of MCMC dates back to the 50's with the article of Metropolis, Rosenbluth, Rosenbluth, Teller and Teller (1953), but it is only with today's computing power that these methods can give their full measure. Like in the classical Monte Carlo schemes, one way of getting insight about the above convergence is by looking at the first moment $\mathbb{E}[S_n]$ and the (asymptotic) variance $\lim n^{-2} \mathbb{V}[S_n]$. There is abundant literature covering these matters—see, for example, Peskun (1973) and Smith and Roberts (1993). A related problem is also to study the rate at which the chain approaches stationarity. Instead, our concern will be the stationary large deviation probabilities

$$\mathbb{P}_\pi[S_n \geq n(\mu + \varepsilon)]. \tag{1}$$

Received May 2001; revised June 2003.
[1]Supported in part by FONDECYT project 1990345 and Cátedra Presidencial de Sistemas Cuánticos 1999.
[2]Supported in part by NSERC.

As $\mu$ and $f$ are arbitrary, this also covers the case $S_n \leq (\mu - \varepsilon)$ so we can restrict ourselves to upper deviations without loss of generality. It is a well-known result from large deviation theory [see Dembo and Zeitouni (1998)] that the asymptotic rate of convergence to zero in (1) is exponential with rate function $I_\eta(x) = \sup_{t \in \mathbb{R}} \{tx - \log \eta(t)\}$, where $\eta(t)$ is the largest eigenvalue of the matrix with coefficients $P(i,j)e^{tf(j)}$. On the other hand, the literature dealing with fixed sample size upper bounds for the above probability is scarce [see Gillman (1993), Dinwoodie (1995) and Lézaud (1998)] and the results do not compare well with the classical bounds when restricted to the independent case [see Hoeffding (1963)].

The above authors use perturbation theory for linear operators to estimate the Perron–Frobenius eigenvalue $\eta(t)$ and obtain upper bounds from the Markov inequality through the matrix representation of the moment generating function $\mathbb{E}[\exp(tS_n)]$. In particular, Lézaud (1998) obtains results for nonreversible and continuous chains.

Our approach contains some elements of the later, but achieves to reduce the initial problem to a simpler one where exact calculations can be carried out. Our bounds are optimal in the sense that the exponential rate is reached asymptotically for a class of Markov chains.

THEOREM 1. *For all pairs $((X_n), f)$, such that $(X_n)$ is a finite, ergodic and reversible Markov chain in stationary state with second largest eigenvalue $\lambda$ and $f$ is a function taking values in $[0,1]$ such that $\mathbb{E}[f(X_i)] = \mu$, the following bounds, with $\lambda_0 = \max(0, \lambda)$, hold for all $\varepsilon > 0$ such that $\mu + \varepsilon < 1$ and all time $n$*

$$\mathbb{P}_\pi[S_n \geq n(\mu + \varepsilon)]$$

$$(2) \qquad \leq \left[\frac{\mu + \bar{\mu}\lambda_0}{1 - 2(\bar{\mu} - \varepsilon)/(1 + \sqrt{\Delta})}\right]^{n(\mu + \varepsilon)} \left[\frac{\bar{\mu} + \mu\lambda_0}{1 - 2(\mu + \varepsilon)/(1 + \sqrt{\Delta})}\right]^{n(\bar{\mu} - \varepsilon)}$$

$$(3) \qquad \leq \exp\left\{-2\frac{1 - \lambda_0}{1 + \lambda_0} n\varepsilon^2\right\},$$

*where*

$$\Delta = 1 + \frac{4\lambda_0(\mu + \varepsilon)(\bar{\mu} - \varepsilon)}{\mu\bar{\mu}(1 - \lambda_0)^2}, \qquad \bar{\mu} = 1 - \mu.$$

In particular, the upper bound given by expression (3) is the large deviation rate function for a two-state chain, which, for $\lambda = 0$, coincides with Hoeffding's bound. The bounds are optimal for $\lambda \geq 0$.

The paper is organized as follows. We will first solve the two-state case, which turns out to be the extremal case needed in the sequel. Next, we handle the case where the cardinality of $E$ is finite by introducing a modification of the spectre of the transition matrix resulting in a new chain $(\tilde{X}_k)$, which will serve as a bridge between the initial chain $(X_k)$ and the two-state case through a positive semidefinitness argument and a convex majorization result. We also compare our bound with existing bounds. Finally, countable chains with a spectral gap can be handled in the same manner.

**1. Solution of the two-state case.** Let $(Y_k)$ be an ergodic Markov chain with state space $\{0,1\}$, transition matrix with second largest eigenvalue $\lambda$ and stationary distribution $\boldsymbol{\mu} = (\bar{\mu}, \mu)'$, $\bar{\mu} + \mu = 1$, $0 < \mu < 1$. Let I be the identity matrix and $\mathbf{1} = (1,1)'$. It is easily seen that $\lambda$ and $\boldsymbol{\mu}$ completely specify the transition mechanism, which is then given by

$$M(\lambda, \boldsymbol{\mu}) = \lambda \mathrm{I} + (1-\lambda)\mathbf{1}\boldsymbol{\mu}'.$$

Following a classical recipe, we derive a bound for the upper deviation probabilities of the empirical sums $S_n = \sum_{k=1}^n Y_k$ using Markov's inequality: for all $t \geq 0$,

(4) $$\mathbb{P}_\mu[S_n \geq n(\mu+\varepsilon)] \leq e^{-nt(\mu+\varepsilon)} \mathbb{E}_\mu[\exp(tS_n)].$$

The expectation on the right-hand side admits the representation [see Dinwoodie (1995)]

(5) $$\boldsymbol{\mu}'[M(\lambda,\boldsymbol{\mu})D_t^2]^n \mathbf{1},$$

where $D_t = \mathrm{diag}(1, e^{t/2})$. Let $D = \mathrm{diag}(\sqrt{\bar{\mu}}, \sqrt{\mu})$ and denote by $\Gamma$ the orthogonal matrix with columns $\boldsymbol{\gamma}_1 = (\sqrt{\bar{\mu}}, \sqrt{\mu})'$ and $\boldsymbol{\gamma}_2 = (-\sqrt{\mu}, \sqrt{\bar{\mu}})'$. The expression (5) then admits the following symmetric form:

$$\boldsymbol{\gamma}_1' D_t G_t^{n-1} D_t \boldsymbol{\gamma}_1,$$

where $G_t = D_t[\boldsymbol{\gamma}_1\boldsymbol{\gamma}_1' + \lambda(\mathrm{I} - \boldsymbol{\gamma}_1\boldsymbol{\gamma}_1')]D_t$. All the above expressions are derived from (5) using the spectral representation $M(\lambda, \boldsymbol{\mu}) = D^{-1}\Gamma \mathrm{diag}(1,\lambda)\Gamma' D$. We will perform this derivation in complete detail for the general case in the next section. Now, the largest eigenvalue $\theta(t)$ of $G_t$ satisfies $\theta(t)^k = \sup_{\|\mathbf{x}\|=1} \|G_t^k \mathbf{x}\|$ and an application of Cauchy–Schwarz's inequality on the last display yields

(6) $$\begin{aligned}\boldsymbol{\gamma}_1' D_t G_t^{n-1} D_t \boldsymbol{\gamma}_1 &\leq \|D_t \boldsymbol{\gamma}_1\| \|G_t^{n-1} D_t \boldsymbol{\gamma}_1\| \\ &\leq \|D_t \boldsymbol{\gamma}_1\|^2 \theta(t)^{n-1}.\end{aligned}$$

PROPOSITION 1. *For the two-state Markov chain with transitions $M$, the stationary upper-deviation probabilities satisfy*

(7) $$\mathbb{P}_\mu[S_n \geq n(\mu+\varepsilon)] \leq \|D_t \boldsymbol{\gamma}_1\|^2 \theta(t)^{-1} \exp\{-n[t(\mu+\varepsilon) - \log\theta(t)]\}.$$

*When $\lambda \geq 0$, it can further be shown that*

$$\mathbb{P}_\mu[S_n \geq n(\mu + \varepsilon)]$$

$$(8) \qquad \leq \left[\frac{\mu + \bar{\mu}\lambda}{1 - 2(\bar{\mu} - \varepsilon)/(\sqrt{\Delta} + 1)}\right]^{n(\mu+\varepsilon)} \left[\frac{\bar{\mu} + \mu\lambda}{1 - 2(\mu + \varepsilon)/(\sqrt{\Delta} + 1)}\right]^{n(\bar{\mu}-\varepsilon)}$$

$$(9) \qquad \leq \exp\left[-2\frac{1-\lambda}{1+\lambda}n\varepsilon^2\right],$$

*where*

$$\Delta = 1 + \frac{4\lambda(\mu + \varepsilon)(\bar{\mu} - \varepsilon)}{\mu\bar{\mu}(1-\lambda)^2}.$$

PROOF. Inequality (7) is obtained from (4) and (6) after rewriting the exponential part. Under the additional condition $\lambda \geq 0$, the nonexponential term in (7) is less than 1. Indeed,

$$\theta(t) \geq \frac{\boldsymbol{\gamma}_1' D_t}{\|D_t\boldsymbol{\gamma}_1\|} G_t \frac{D_t\boldsymbol{\gamma}_1}{\|D_t\boldsymbol{\gamma}_1\|}$$

$$= \|D_t\boldsymbol{\gamma}_1\|^2 + \lambda\frac{\boldsymbol{\gamma}_1' D_t^2(\,\mathrm{I} - \boldsymbol{\gamma}_1\boldsymbol{\gamma}_1')D_t^2\boldsymbol{\gamma}_1}{\|D_t\boldsymbol{\gamma}_1\|^2}$$

$$\geq \|D_t\boldsymbol{\gamma}_1\|^2,$$

since $\mathrm{I} - \boldsymbol{\gamma}_1\boldsymbol{\gamma}_1'$ is positive semidefinite. Then we have

$$\mathbb{P}_\mu[S_n \geq n(\mu + \varepsilon)] \leq \exp\{-n[t(\mu + \varepsilon) - \log\theta(t)]\}.$$

Taking the infimum for $t \geq 0$, we obtain

$$\mathbb{P}_\mu[S_n \geq n(\mu + \varepsilon)] \leq \exp\{-nI_\theta(\mu + \varepsilon)\},$$

where $I_\theta(x) = \sup_{t \in \mathbb{R}}\{tx - \log\theta(t)\}$ [see Dembo and Zeitouni (1998), Lemma 2.2.5] is the rate function of the empirical averages $n^{-1}S_n$. An explicit computation of this rate function will prove (8). A simple calculation yields

$$(10) \qquad \theta(t) = \tfrac{1}{2}[\mathrm{Tr}(G_t) + \sqrt{\mathrm{Tr}^2(G_t) - 4\lambda e^t}],$$

where $\mathrm{Tr}(\cdot)$ denotes the trace. Taking $0 < x < 1$ arbitrary but fixed and looking for the zeros of the equation $\frac{\partial}{\partial t}[tx - \log\theta(t)] = 0$, we obtain

$$(11) \qquad (x - \bar{x})\sqrt{\mathrm{Tr}^2(G_t) - 4\lambda e^t} - [(\mu + \bar{\mu}\lambda)e^t - (\bar{\mu} + \mu\lambda)] = 0.$$

Multiplying by the conjugate, simplifying and expanding into powers of $e^t$, a quadratic polynomial emerges whose roots are possible candidates for the maximum value $t_0$, which is found to be the following (see Appendix A for the details):

$$t_0 = \log\left[\frac{(\bar{\mu} + \mu\lambda)[\sqrt{\Delta} - (\bar{x} - x)]}{(\mu + \bar{\mu}\lambda)[\sqrt{\Delta} + (\bar{x} - x)]}\right],$$

where

$$\Delta = \left[1 + \frac{4\lambda x\bar{x}}{\mu\bar{\mu}(1-\lambda)^2}\right] \quad \text{and} \quad \bar{x} = 1 - x.$$

Now we can evaluate expression (10) to determine the Perron eigenvalue $\theta(t_0)$, which we call simply $\theta$,

$$\theta = \frac{(\bar{\mu} + \mu\lambda)[\sqrt{\Delta} + 1]}{\sqrt{\Delta} + \bar{x} - x}.$$

This yields the rate function $I_\theta$, which after simplifying can be written as

$$(12) \quad I_\theta(x) = -x\log\left[\frac{\mu + \bar{\mu}\lambda}{1 - 2\bar{x}/(\sqrt{\Delta}+1)}\right] - \bar{x}\log\left[\frac{\bar{\mu} + \mu\lambda}{1 - 2x/(\sqrt{\Delta}+1)}\right]$$

and the right-hand side of (8) is just the explicit form of $\exp\{-nI_\theta(\mu + \varepsilon)\}$. To prove the uniform bound (9), we will show that $I_\theta(x) \geq \frac{(1-\lambda)}{(1+\lambda)}(x - \mu)^2$. First, we differentiate $g(x) = I_\theta(x)/(x-\mu)^2$ to obtain $g'(x) = (x-\mu)^{-3}h(x)$, where $h(x) = (x-\mu)I'_\theta(x) - 2I_\theta(x)$. Studying the two first derivatives of the numerator $h(x)$, it can be shown (see Appendix B) that

$$h(x) \begin{cases} > 0, & \text{if } |x - 1/2| > |\mu - 1/2|, \\ = 0, & \text{if } x = \mu \text{ or } x = \bar{\mu}, \\ < 0, & \text{if } |x - 1/2| < |\mu - 1/2|, \end{cases}$$

and $g'$ has the following behavior:

$$g'(x) \begin{cases} < 0, & \text{if } x < \bar{\mu}, \\ = 0, & \text{if } x = \bar{\mu}, \\ > 0, & \text{if } x > \bar{\mu}, \end{cases}$$

from what we deduce that $g$ attains a global minimum at $x = \bar{\mu}$. Now, a Taylor expansion for $g(\bar{\mu}) = I_\theta(\bar{\mu})(\bar{\mu} - \mu)^{-2}$ in terms of $r = (\bar{\mu} - \mu)(1 - \lambda)(1+\lambda)^{-1}$ gives

$$g(\bar{\mu}) = (\bar{\mu} - \mu)^{-1}\log[(1+r)/(1-r)]$$

$$= 2\frac{1-\lambda}{1+\lambda}r^{-1}\left(r + \frac{1}{3}r^3 + \frac{1}{5}r^5 + \cdots\right)$$

$$\geq 2\frac{1-\lambda}{1+\lambda},$$

and from the definition of $g$ we then have

$$I_\theta(x) \geq g(\bar{\mu})(x-\mu)^2 \geq 2\frac{1-\lambda}{1+\lambda}(x-\mu)^2. \qquad \square$$

**2. General case.** Let $(X_k)$ have transitions $P = (p_{ij})$ with stationary distribution $\boldsymbol{\pi} = (\pi_i)$ such that

$$\pi_i p_{ij} = \pi_j p_{ji} \tag{13}$$

for all $i, j$ in the finite and ordered state space $(E, \leq)$. (It is convenient to leave the order $\leq$ unspecified.) Consider now a bounded function $f : E \to \mathbb{R}$ with $\min f(E) = 0, \max f(E) = 1$ and such that the stationary mean $\mathbb{E}_\pi[f(X_k)]$ is equal to a fixed number $\mu$. Applying an obvious affine transformation we can always set $\min f(E) = 0$, and $\max f(E) = 1$; this has no bearing on our argumentation and will only make the expressions more concise.

Using condition (13) we now derive a spectral decomposition for the matrix $P$. This result is well known [see, e.g., Green and Han (1992)], but our derivation contains some new elements and will allow us to introduce most of the notation needed for the sequel in a smooth way. Call $D$ the diagonal matrix with (diagonal) elements $\sqrt{\pi_i}$, where $i$ runs through $(E, \leq)$. Condition (13) says that $D^2 P$ is symmetric, and hence, $DPD^{-1}$ is symmetric too. Since the last product shares its spectrum with $P$, the transition matrix has real eigenvalues $\lambda_l$, and further admits the spectral representation

$$P = D^{-1} \Gamma \operatorname{diag}[\lambda_l] \Gamma^t D, \tag{14}$$

where $\Gamma$ is orthogonal. Furthermore, since $P$ is irreducible and aperiodic, from the Perron–Frobenius theorem we know that the largest eigenvalue $\lambda_1 = 1$ strictly dominates in modulus any other eigenvalue and also the corresponding eigenvector in (14) is positive. In fact, using (14) and stationarity we get

$$D^{-1} \Gamma \operatorname{diag}[\delta_1^l] \Gamma^t D = A,$$

where $A$ is the limiting matrix $(A_{ij}) = \pi_j$ and $\delta_1^l$ is Kronecker's delta. So far, all this is well known. Now, let $\lambda$ be the second largest eigenvalue of $P$, $\lambda_0 = \max(0, \lambda)$ and consider

$$\begin{aligned} Q &= D^{-1} \Gamma \operatorname{diag}[\max(\lambda_0, \lambda_l)] \Gamma^t D \\ &= \lambda_0 I + (1 - \lambda_0) A. \end{aligned} \tag{15}$$

Clearly the rows of $Q$ sum to 1 and the off-diagonal elements are positive, hence it is a good candidate to be a stochastic matrix, with $Q$ being a convex linear combination of stochastic matrices. It only remains to show that the diagonal elements are always nonnegative. Since $D(Q - P)D^{-1}$ is positive semidefinite, the desired result follows from

$$q_{ii} - p_{ii} = \mathbf{e}_i^t D(Q - P) D^{-1} \mathbf{e}_i \geq 0,$$

where $\mathbf{e}_i$ is the corresponding canonical vector. Observe that since the map $\lambda_l \mapsto \max(\lambda_0, \lambda_l)$ leaves the largest eigenvalue unchanged, $\pi$ is stationary for $Q$-alternatively check it directly from (15). When $\lambda \geq 0$ we have $\lambda = \lambda_0$ so the largest and the second largest eigenvalues of $Q$ and the ones of $P$ are the same. Since $\text{Tr}(Q) = 1 + \lambda(1 - |E|) \geq 0$, where $|E|$ is the cardinality of $E$, we have that $\lambda \geq -(|E|-1)^{-1}$; hence when dealing with arbitrarily large chains, $\lambda_0$ and $\lambda$ cannot be very far apart. A crucial property of these transitions is the preservation of the Markov property under any transformation.

Here we give a simple construction for deriving a Markov chain $(X'_k)$ with transition probabilities $Q$. Consider $\lambda_0 \in [0,1)$, $\pi$ and $E$ as fixed but arbitrary and let $(I_k)$ and $(Z_k)$ be independent sequences of i.i.d. random variables with respective distribution Bernoulli$(1 - \lambda_0)$ and $\pi$ on $E$. Let $I'_1 = 1$ and $I'_k = I_k$ for $k > 1$. It is easy to verify that the construction

$$X'_k = \sum_{\{j:\, 1 \leq j \leq k\}} \left\{ \prod_{\{\ell:\, j < \ell \leq k\}} (1 - I'_\ell) \right\} I'_j Z_j$$

works. Moreover, if we set

$$N(j) = \sum_{\{k:\, j \leq k \leq n\}} \left\{ \prod_{\{\ell:\, j < \ell \leq k\}} (1 - I'_\ell) \right\} I'_j,$$

then

(16) $$\sum_{k=1}^n f(X'_k) = \sum_{j=1}^n N(j) f(Z_j)$$

with independence between $N(1), \ldots, N(n)$ and $f(Z_1), \ldots, f(Z_n)$. It is plain from this representation that applying a transformation on the observations only amounts to changing the distribution of the i.i.d. $Z$'s, that is, changing $E$ and $\pi$. Since $(X'_k)$ is Markovian regardless of $E$ and $\pi$, any transformation will preserve the Markov property.

Our goal now will be to relate the moment generating function $\mathbb{E}_\pi[tS_n]$ to its 0–1 counterpart studied previously. This will be done in two steps. First, we compare $S_n = \sum_{k=1}^n f(X_k)$ with the $n$th partial sum $S'_n$ of the chain $(f(X'_k))$, where $(X'_k)$ has transitions $Q$. Second, we establish a stochastic majorization property for $(X'_k)$ that will enable us to relate it to the two-state case.

2.1. *Step* 1. As seen previously for the two-state case, the moment-generating function of the partial sums $S_n$ can be written as

$$\mathbb{E}_\pi[\exp(tS_n)] = \boldsymbol{\pi}'[PD_t^2]^n \mathbf{1},$$

where $D_t$ is the diagonal matrix with entries $\exp(tf(i)/2)$. Since $\pi$ is stationary for $P$ and since diagonal matrices commute, from the spectral representation (14) we get

$$\boldsymbol{\pi}'[PD_t^2]^n \mathbf{1} = \boldsymbol{\pi}' D_t [D_t P D_t]^{n-1} D_t \mathbf{1}$$

$$= \boldsymbol{\pi}' D_t D^{-1} [D_t \Gamma \operatorname{diag}[\lambda_l] \Gamma^t D_t]^{n-1} D D_t \mathbf{1}$$
$$= \boldsymbol{\gamma}_1' D_t G_t^{n-1} D_t \boldsymbol{\gamma}_1,$$

where $\boldsymbol{\gamma}_1 = (\sqrt{\pi_i})$ is the first column of $\Gamma$ and $G_t = D_t \Gamma \operatorname{diag}[\lambda_l] \Gamma^t D_t$. Since $G_t$ has nonnegative entries and is irreducible, from the Perron–Frobenius theorem the largest eigenvalue $\zeta(t)$ satisfies $\zeta(t)^k = \sup_{\|x\|=1} \|G_t^k\|$ and the same argument that we used to obtain (7) yields

$$\boldsymbol{\gamma}_1' D_t G_t^{n-1} D_t \boldsymbol{\gamma}_1 \leq \|D_t \boldsymbol{\gamma}_1\|^2 \zeta(t)^{n-1}.$$

If we introduce $H_t = D_t \Gamma \operatorname{diag}[\max(\lambda_0, \lambda_i)] \Gamma^t D_t$, and denote its largest eigenvalue by $\eta(t)$, then $H_t - G_t$ is positive semidefinite and $\eta(t)$ dominates $\zeta(t)$ so that

$$\|D_t \boldsymbol{\gamma}_1\|^2 \zeta(t)^{n-1} \leq \|D_t \boldsymbol{\gamma}_1\|^2 \eta(t)^{n-1}.$$

PROPOSITION 2. *The large deviation probabilities satisfy*

(17) $\quad \mathbb{P}_\mu[S_n \geq n(\mu + \varepsilon)] \leq \|D_t \boldsymbol{\gamma}_1\|^2 \zeta(t)^{-1} \exp\{-n[t(\mu + \varepsilon) - \log \zeta(t)]\}$

(18) $\quad\quad\quad\quad\quad\quad\quad \leq \|D_t \boldsymbol{\gamma}_1\|^2 \eta(t)^{-1} \exp\{-n[t(\mu + \varepsilon) - \log \eta(t)]\}.$

*Since $\lambda_0 \geq 0$, we further have*

$$\mathbb{P}_\mu[S_n \geq n(\mu + \varepsilon)] \leq \exp\{-n[t(\mu + \varepsilon) - \log \eta(t)]\}.$$

PROOF. Markov's inequality together with (17) and (18) imply the two first inequalities. Now, just as in the 0–1 case, the condition $\lambda_0 \geq 0$ guarantees $\eta(t) \geq \|D_t \boldsymbol{\gamma}_1\|^2$ and the last inequality holds. □

2.2. *Step* 2.

THEOREM 2. *Let $(X_k')$ have transitions $Q$ and let $\mathbb{E}_\pi[f(X_k')] = \mu$. Then for any convex function $\Psi : \mathbb{R} \to \mathbb{R}$, we have*

$$\mathbb{E}_\pi[\Psi(f(X_1') + \cdots + f(X_n'))] \leq \mathbb{E}_\mu[\Psi(Y_1 + \cdots + Y_n)],$$

*where $(Y_k)$ is the two-state chain with transitions $M(\lambda_0, \boldsymbol{\mu})$.*

PROOF. The proof is based on the representation (16) and a construction. We introduce random variables $B_j$, $j = 1, \ldots, n$ and we consider a joint distribution on $(B_j, Z_j)$ such that the conditional distribution of $B_j$, given $Z_j$, is a Bernoulli distribution with mean $f(Z_j)$. We assume that $(B_1, Z_1), \ldots, (B_n, Z_n)$ are independent. Therefore, the marginal distribution of $B_j$ is Bernoulli($\mu$) for $j = 1, \ldots, n$. As in expression (16) we set

$$Y_k = \sum_{\{j\,:\,1 \leq j \leq k\}} \left\{ \prod_{\{\ell\,:\,j < \ell \leq k\}} (1 - I'_\ell) \right\} I'_j B_j,$$

so

$$\sum_{k=1}^{n} Y_k = \sum_{j=1}^{n} N(j) B_j.$$

Jensen's inequality combined with a conditional expectation tells us that

$$\Psi\left(\sum_{j=1}^{n} N(j) f(Z_j)\right) \leq \mathbb{E}_\pi\left[\Psi\left(\sum_{j=1}^{n} N(j) B_j\right) \bigg| N, Z\right].$$

Taking expectation on both sides, we obtain that

$$\mathbb{E}_\pi[\Psi(S'_n)] \leq \mathbb{E}_\mu\left[\Psi\left(\sum_{k=1}^{n} Y_k\right)\right]. \qquad \square$$

REMARK 1. In the above proof, it is implicit that the endpoints of $f(E)$ are $a = 0$ and $b = 1$. When $a < b$ are arbitrary, the corresponding extremal chain lives in $\{a, b\}$ and the transitions are determined by

$$m_{ab} = (1 - \lambda_0) \frac{\mu - a}{b - a}, \qquad m_{ba} = (1 - \lambda_0) \frac{b - \mu}{b - a}.$$

Theorem 2 deals with *stochastic ordering* of random variables. The particular stochastic order used here is known in the literature as the convex ordering: $X \preceq Y$ if $\mathbb{E}[\Psi(X)] \leq \mathbb{E}[\Psi(Y)]$ for all convex real valued $\Psi$ (such that the expectations exist). The result can be stated as $X'_1 + \cdots + X'_n \preceq Y_1 + \cdots + Y_n$. Observe that under stationarity $X'_k \preceq Y_k$, so we have transition schemes $Q$ and $M$ under which the stochastic order relation is preserved for the respective partial sums. When there is independency within each sequence, it is known that the convex ordering of the marginals implies the same ordering for the corresponding partial sums [see Marshall and Olkin (1979)]; our result shows that this preservation property can occur in the Markovian setting as well.

We now have all the necessary tools to prove our main result.

PROOF OF THEOREM 1. Applying Theorem 2 with $\Psi(x) = \exp(tx)$, we get

$$\eta(t) \leq \lim_{n \to \infty} n^{-1} \log \mathbb{E}_\pi\{\exp(tS'_n)\}$$

$$\leq \lim_{n\to\infty} n^{-1} \log \mathbb{E}_\mu \left\{ \exp\left(t \sum_{k=1}^n Y_k\right) \right\}$$

$$= \theta(t)$$

and then from Proposition 2 we obtain

$$\mathbb{P}_\pi[S_n \geq n(\mu + \varepsilon)] \leq \inf_{t \geq 0} \exp\{-n[t(\mu+\varepsilon) - \log \eta(t)]\}$$

$$\leq \inf_{t \geq 0} \exp\{-n[t(\mu+\varepsilon) - \log \theta(t)]\}$$

$$= \exp\{-n I_\theta(\mu+\varepsilon)\},$$

where $I_\theta$ is the rate function (12). The stated upper bounds then follow from Proposition 1. □

REMARK 2. With a little more effort we can see that we do not need to assume that $P$ is aperiodic in Theorem 1. Indeed, given periodic but irreducible and reversible $P$, it is possible to construct a sequence of aperiodic, irreducible and reversible chains $P_m$, such that $P_m$ converges to $P$ as $m$ tends to infinity. Since $\lambda$, $\pi$ and $\mu$ are continuous functions of $P$, Theorem 1 will hold for all $P_m$ and the result will hold for $P$ as well, by continuity. In fact, an eigenvalue near of even equal to $-1$ is not a problem as only Césaro sums $S_n$ are considered here.

REMARK 3. Following Remark 1, the theorem remains valid when the end points $a < b$ of $f(E)$ are arbitrary. In this case the values $\mu$ and $\varepsilon$ in the bounds are to be replaced by $\frac{\mu-a}{b-a}$ and $\frac{\varepsilon}{b-a}$, respectively.

It is clear from the proof of Theorem 1 that, under the condition $\lambda \geq 0$, the rate functions $I_\zeta(x) = \sup_{t \in \mathbb{R}}\{tx - \log \zeta(t)\}$, $I_\eta(x) = \sup_{t \in \mathbb{R}}\{tx - \log \eta(t)\}$ and $I_\theta(x) = \sup_{t \in \mathbb{R}}\{tx - \log \theta(t)\}$ corresponding to $S_n, S'_n$ and $\sum_{k=1}^n Y_k$, respectively, satisfy

(19) $$I_\zeta(x) \geq I_\eta(x) \geq I_\theta(x).$$

When $P = Q$, there is equality on the leftmost side. Furthermore, when $f$ is 0–1, we have equality in the rightmost side. Hence, when $P = Q$ and $f : E \to \{0, 1\}$, the exponential rate given in our first upper bound cannot be improved upon. In particular, when the chain is independent, the theorem yields the well-known Hoeffding's inequality

$$\mathbb{P}[S_n \geq n(\mu+\varepsilon)] \leq \left(\frac{\mu}{\mu+\varepsilon}\right)^{n(\mu+\varepsilon)} \left(\frac{\bar\mu}{\bar\mu - \varepsilon}\right)^{n(\bar\mu-\varepsilon)}.$$

REMARK 4. A closer look reveals that the leftmost inequality in (19) is true for all $\lambda$. This suggests the possibility that the theorem might be true for all admissible values of $\lambda$. But this is not so, numerical evidence show that the bounds do not hold without the condition $\lambda \geq 0$.

**3. Comparisons.** Gillman (1993) was the first to obtain a finite sample size exponential bound for the large deviation probabilities using perturbation theory. Successive refinements of the technique allowed Dinwoodie (1995) and Lézaud (1998) to improve this bound. Among these, the later work contains the best results and we shall use them for the comparisons. With $f$ satisfying our usual assumptions, Theorem 1.1 of Lézaud (1998) gives in our particular case

$$(20) \qquad \mathbb{P}_\pi[S_n \geq n(\mu + \varepsilon)] \leq e^{(1-\lambda)/5} \exp\left\{-\frac{(1-\lambda)n\varepsilon^2}{4\mu[1 + h(5\varepsilon/\mu)]}\right\},$$

where $h(x) = \sqrt{1-x} - (1-x)$. Let us denote $L(\mu, \varepsilon)$ the exponential rate in (20) and $I_\theta(\mu + \varepsilon)$ is the exponential rate in the bound (2). Observe first that since $I_\theta(\mu + \varepsilon)$ comes from the rate function of the two-state case, and since $[\mu + \varepsilon, \infty)$ is a continuity set of $I_\theta$, then sampling from this chain implies

$$\lim_{n \to \infty} n^{-1} \log \mathbb{P}_\mu[S_n \geq n(\mu + \varepsilon)] = -I_\theta(\mu + \varepsilon) \leq -L(\mu, \varepsilon),$$

hence, when $\lambda \geq 0$, the rate $I_\theta(\mu + \varepsilon)$ always yields a better bound. A limited Taylor expansion of $I_\theta$ around $\mu$ gives an idea of the ratio of these quantities

$$\frac{I_\theta(\mu + \varepsilon)}{L(\mu, \varepsilon)} = \frac{2}{\bar{\mu}(1+\lambda)} + o(\varepsilon).$$

## APPENDIX A

The leading, middle and constant terms of the convex quadratic polynomial obtained from (15) are

$$a = [1 - (2x-1)^2](\mu + \bar{\mu}\lambda),$$
$$b = -2\{[\mu\bar{\mu}(1-\lambda)^2 + \lambda][1 + (2x-1)^2] - 2\lambda(2x-1)^2\}$$

and

$$c = (\bar{\mu} + \mu\lambda)^2[1 - (2x-1)^2],$$

respectively. After some simplifications the discriminant $b^2 - 4ac$ can be written as

$$16(2x-1)^2[\mu\bar{\mu}(1-\lambda)^2]^2\left[1 + \frac{4\lambda x(1-x)}{\mu\bar{\mu}(1-\lambda)^2}\right] \geq 0$$

and the roots $t_0^+, t_0^-$ are given, respectively, by

$$
\begin{aligned}
(21) \quad & \frac{\mu\bar{\mu}(1-\lambda)^2[1+(2x-1)^2]+\lambda[1-(2x-1)^2]}{[1-(2x-1)^2](\mu+\bar{\mu}\lambda)} \\
& \pm \frac{2(2x-1)\mu\bar{\mu}(1-\lambda)^2\sqrt{\Delta}}{[1-(2x-1)^2](\mu+\bar{\mu}\lambda)}.
\end{aligned}
$$

Now, consider the conjugate product

$$
\begin{aligned}
& [\sqrt{\Delta}+(1-2x)][\sqrt{\Delta}-(1-2x)] \\
& = \frac{(\mu+\bar{\mu}\lambda)(\bar{\mu}+\mu\lambda)[1-(2x-1)^2]}{\mu\bar{\mu}(1-\lambda)^2}.
\end{aligned}
$$

Since it is positive for $0 < x < 1$ and since both terms on the left-hand side are positive at $x = 1/2$ and continuous, each is positive for all $0 < x < 1$. Multiplying the numerator and denominator in the expression (21) by $\sqrt{\Delta}+(1-2x)$, for $t_0^+$, and by $\sqrt{\Delta}-(1-2x)$, for $t_0^-$, the roots can be written as

$$
\frac{(\bar{\mu}+\mu\lambda)[\sqrt{\Delta}-(1-2x)]}{(\mu+\bar{\mu}\lambda)[\sqrt{\Delta}+(1-2x)]}, \qquad \frac{(\bar{\mu}+\mu\lambda)[\sqrt{\Delta}+(1-2x)]}{(\mu+\bar{\mu}\lambda)[\sqrt{\Delta}-(1-2x)]}.
$$

Except for $x = 1/2$, where they coincide, exactly one of these is the solution of (11), the other being the solution to the conjugate equation. To arbitrate, let us, evaluate the rightmost term in (11) for the first candidate. We obtain

$$
(\mu+\bar{\mu}\lambda)\frac{(\bar{\mu}+\mu\lambda)[\sqrt{\Delta}-(1-2x)]}{(\mu+\bar{\mu}\lambda)[\sqrt{\Delta}+(1-2x)]} - (\bar{\mu}+\mu\lambda) = \frac{2(2x-1)(\bar{\mu}+\mu\lambda)}{\sqrt{\Delta}+(1-2x)}.
$$

Since this expression shares its sign with the leftmost term in (11), we have found the maximizing value.

## APPENDIX B

The behavior of $h(x) = (x-\mu)I_\theta'(x) - 2I_\theta(x)$ depends on whether $\mu < 1/2$ or $\mu \geq 1/2$; we shall carry out the analysis for the first case, the other being similar but somewhat less involved. To begin with, the first derivatives of $I_\theta$ are found to be

$$
\begin{aligned}
I_\theta'(x) &= -\log\left[\frac{\mu+\bar{\mu}\lambda}{1-2\bar{x}/(1+\sqrt{\Delta})}\right] + \log\left[\frac{\bar{\mu}+\mu\lambda}{1-2x/(1+\sqrt{\Delta})}\right], \\
I_\theta''(x) &= (\sqrt{\Delta}x\bar{x})^{-1}, \\
I_\theta^{(3)}(x) &= \frac{(x-\bar{x})(3\Delta-1)}{2\Delta^{3/2}(x\bar{x})^2}
\end{aligned}
$$

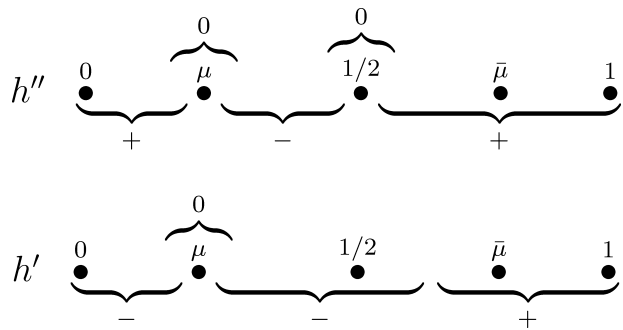

Fig. 1.

so that $I_\theta(\mu) = I'_\theta(\mu) = 0$, while $I_\theta^{(3)}(x) \propto (x - \bar{x})$, since the other terms are positive. Now, we have $h'(x) = (x - \mu)I''_\theta(x) - I'_\theta(x)$, $h''(x) = (x - \mu)I_\theta^{(3)}(x)$, and Figure 1 summarizes the analysis of their sign.

Combining this with the fact that $h(\mu) = h(\bar{\mu}) = 0$, we see that these are the only zeros and further, $h$ is negative in $(\mu, \bar{\mu})$ and positive in $(0, \mu) \cup (\bar{\mu}, 1)$.

**Acknowledgment.** We would like to thank the anonymous referee whose helpful comments allowed us to improve the presentation of the paper.

Depto. de Ingeniería Matemática
Universidad de Concepción
Avda. Esteban S. Iturra s/n
Casilla 160-C, Concepción
Chile
e-mail: caleon@mat.ulaval.ca

Dépt. de Mathématiques et
  Statistique
Université de Montréal
C.P. 6128, succ. Centre Ville
Montréal, Quebec H3C 3J7
Canada
e-mail: perronf@dms.umontreal.ca